\documentclass[12pt]{amsart}

\renewcommand{\b}{\mathbf{b}}

\newcommand{\om}{\omega}

\DeclareFontFamily{U}{mathx}{\hyphenchar\font45}
\DeclareFontShape{U}{mathx}{m}{n}{
      <5> <6> <7> <8> <9> <10>
      <10.95> <12> <14.4> <17.28> <20.74> <24.88>
      mathx10
      }{}
\DeclareSymbolFont{mathx}{U}{mathx}{m}{n}
\DeclareFontSubstitution{U}{mathx}{m}{n}
\DeclareMathAccent{\widecheck}{0}{mathx}{"71}

\renewcommand{\kill}[1]{}
\newcommand{\dummy}[1]{\mbox{}}

\makeatletter
\newcommand{\xequal}[2][]{\ext@arrow 0055{\equalfill@}{#1}{#2}}
\def\equalfill@{\arrowfill@\Relbar\Relbar\Relbar}
\makeatother

\newcommand{\1}{\ensuremath{\ol{\mathrm{P}}}}

\newcommand{\bca}{\begin{cases}}
\newcommand{\eca}{\end{cases}}

\newcommand{\bdet}{\te{bdet}}

\newcommand{\bpic}{\begin{picture}}\newcommand{\epic}{\end{picture}}

\newcommand{\beda}{\begin{edaenumerate}}
\newcommand{\eeda}{\end{edaenumerate}}

\newcommand{\q}{\quad}

\newcommand{\bq}{\begin{quote}}\newcommand{\eq}{\end{quote}}

\newcommand{\be}{\begin{enumerate}}\newcommand{\ee}{\end{enumerate}}
\newcommand{\bce}{\begin{center}}\newcommand{\ece}{\end{center}}
\newcommand{\bde}{\begin{description}}\newcommand{\ede}{\end{description}}
\newcommand{\bri}{\begin{flushright}}\newcommand{\eri}{\end{flushright}}
\newcommand{\bb}{\begin{block}}\newcommand{\eb}{\end{block}}
\newcommand{\bt}{\begin{thm}}\newcommand{\et}{\end{thm}}
\newcommand{\bpf}{\begin{proof}}\newcommand{\epf}{\end{proof}}
\newcommand{\bex}{\begin{ex}}\newcommand{\eex}{\end{ex}}
\newcommand{\bexr}{\begin{exr}}\newcommand{\eexr}{\end{exr}}
\newcommand{\bft}{\begin{fact}}\newcommand{\eft}{\end{fact}}
\newcommand{\brk}{\begin{rmk}}\newcommand{\erk}{\end{rmk}}
\newcommand{\ba}{\begin{align*}}\newcommand{\ea}{\end{align*}}
\newcommand{\bexe}{\begin{exe}}\newcommand{\eexe}{\end{exe}}

\newcommand{\bit}{\begin{itemize}}\newcommand{\eit}{\end{itemize}}

\newcommand{\bcm}{}
\newcommand{\ol}{\overline}

\newcommand{\fr}{\frac}
\newcommand{\cc}{\ensuremath{\mathbf{C}}}

\newcommand{\zz}{\ensuremath{\mathbf{Z}}}

\newcommand{\bd}{\begin{defn}}\newcommand{\ed}{\end{defn}}
\newcommand{\bp}{\begin{prop}}\newcommand{\ep}{\end{prop}}

\newcommand{\sub}{\subseteq}
\newcommand{\lam}{\lambda}

\newcommand{\te}{\text}

\newcommand{\di}{\displaystyle}\renewcommand{\a}{\ensuremath{\bm{a}}}
\renewcommand{\b}{\ensuremath{\bm{b}}}

\newcommand{\x}{\ensuremath{\bm{x}}}

\renewcommand{\b}{\beta}

\renewcommand{\a}{\alpha}

\renewcommand{\int}{\in T}

\renewcommand{\x}{\mathbf{x}}

\usepackage[dvipdfmx]{graphicx}
\usepackage[dvipsnames]{xcolor}
\usepackage{float,afterpage}
\usepackage{asymptote,layout}
\usepackage{wrapfig,epic}

\usepackage{geometry}
\usepackage{exscale,latexsym,bm}
\usepackage{amssymb,enumerate,amsmath,amsthm,amsfonts}
\usepackage{verbatim,fancybox,wasysym,fancyhdr,type1cm}
\usepackage[frame,all,poly,curve,knot,arrow]{xy}
\usepackage{colortbl}
\usepackage{boxedminipage}
\usepackage{multirow}

\graphicspath{{./figs/}}
\theoremstyle{definition}
\newtheorem{thm}{Theorem}[section]
\newtheorem{lem}[thm]{Lemma}
\newtheorem{prop}[thm]{Proposition}\newtheorem{cor}[thm]{Corollary}

\newtheorem{exr}[thm]{Exercise}

\newtheorem{prob}[thm]{Problem}
\newtheorem*{idea}{Main idea}

\newtheorem{ex}[thm]{Example}

\newtheorem{defn}[thm]{Definition}\newtheorem{rmk}[thm]{Remark}
\newtheorem{fact}[thm]{Fact}
\newtheorem{block}[thm]{}
\newtheorem*{exe}{Exercise}

\begin{document}
\title[Determinants and Signed bigrassmannian polynomials] 
{Enumerative combinatorics on determinants and signed bigrassmannian polynomials} 

\author[M. Kobayashi] 
       {Masato Kobayashi}
       \address{Department of Engineering\\
       Kanagawa University, 3-27-1 Rokkaku-bashi, Yokohama 221-8686, Japan.}
\email{kobayashi@math.titech.ac.jp}






\thanks{This article is published as 
Math. J. Okayama Univ. \textbf{57} (2015), 159--172.
}


\begin{abstract}
As an application of linear algebra for enumerative combinatorics, we introduce two new ideas, signed bigrassmannian polynomials and bigrassmannian determinant. First, a signed bigrassmannian polynomial is 
a variant of the statistic given by the number of bigrassmannian permutations below a permutation in Bruhat order as Reading suggested  (2002) and afterward the author developed (2011).
Second, bigrassmannian determinant is a $q$-analog of the determinant with respect to our statistic. It plays a key role for a determinantal expression of those polynomials. We further show that bigrassmannian determinant satisfies weighted condensation as a generalization of Dodgson, Jacobi-Desnanot and Robbins-Rumsey (1986).\end{abstract}
\maketitle

\section{Introduction}

The purpose of this article is to introduce two new ideas, \emph{signed bigrassmannian polynomials} and \emph{bigrassmannian determinant} as an application of linear algebra for enumerative combinatorics. We begin with explaining our motivation.
\subsection{Reading's problem (2002): bigrassmannian statistic}


Permutation statistics has been of great importance in enumerative combinatorics; in particular, Mahonian and Eulerian statistics, such as inversions and descent numbers, are fundamental in the theory.
Here what we deal with is a certain new statistic $\beta$, which we call \emph{bigrassmannian statistic}.
Reading \cite{reading} suggested the following problem:
  
\begin{prob}
Let $\beta(w)$ be the number of join-irreducible (equivalently, bigrassmannian) permutations weakly below a permutation $w$ in Bruhat order. Find its generating function $\sum_{w\in S_n}q^{\beta(w)}$.
\end{prob}
He gave examples of such generating functions for smaller $n$'s:
\[\left.\begin{array}{ccc}
1+q  & \phantom{aaa}  &(S_2)   \\1+2q+2q^3+q^4  &   &  (S_3) \\1+3q+q^2+4q^3+2q^4+2q^5+2q^6+4q^7+q^8+3q^9+q^{10}
  &   &(S_4)  \end{array}\right.\]
Unfortunately, we failed to find any patterns of these coefficients nor factors of such polynomials. Instead, in this article, we study the following \emph{signed} statistic (as signed Mahonian or signed Eulerian statistics):
\[B_n(q)=\sum_{w\in S_n} (-1)^{\ell(w)}q^{\beta(w)}\]
where $\ell(w)$ is the number of inversions; let us call $\{B_n(q)\mid n=1, 2, 3, \dots\}$ \emph{signed bigrassmannian polynomials}.
Fortunately, we could find satisfactory descriptions of such polynomials. It turned out that it is also worthwhile to study these polynomials with a connection to \emph{tournaments} and \emph{Vandermonde determinant}. 
 Since each $B_n(q)$ is a signed sum over the symmetric group, it is natural to come to this idea:





\begin{idea}Use the determinant to find $B_n(q)$.
\end{idea}
The determinant is usually a function which outputs a scalar. For our purpose to find $B_n(q)$, we introduce its $q$-analog (Section \ref{sec4}); we call it \emph{bigrassmannian determinant}.\\
As main results, we will prove three theorems:
\begin{itemize}
\item Theorem \ref{mth1}: a factorization of $B_n(q)$.
\item Theorem \ref{mth15}: a determinantal expression of $B_n(q)$.
\item Theorem \ref{mth2}: weighted condensation for bigrassmannian determinant.
\end{itemize}
In addition, we observe a corollary after each of these theorems.
\subsection{Overview}
In Section \ref{to}, we review some classic results on tournaments and Vandermonde determinant as mentioned above. These facts will play a fundamental role in the sequel.
In Section \ref{too}, we introduce $\beta$-statistic for tournaments as well as permutations. Then we find factors of $B_n(q)$ using weighted Vandermonde determinant. 
Section \ref{sec4} continues to study $B_n(q)$ (from a little different aspect);
we give a definition of bigrassmannian determinant for square matrices as a $q$-analog of the original one. This new idea leads to a determinantal expression of $B_n(q)$ as we shall see. Further, we prove that bigrassmannian determinant satisfies weighted condensation. 
It slightly generalizes the construction of Robbins-Rumsey's $\lam$-determinants \cite{rr}.
We end with some comments for future work in Section \ref{sec5}.
 



\section{Tournaments and Vandermonde determinant}\label{to}
We begin with combinatorics of tournaments and Vandermonde determinant.

\subsection{Tournaments}

\bd{
A \emph{tournament} is a complete digraph with vertices labeled by $1, 2, \dots, n$.
We denote by $T_n$ the set of all tournaments.
}\ed


\bex{Here are eight elements in $T_3$:
\label{t3}
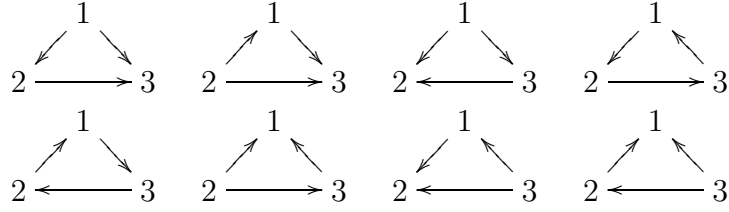
\[\begin{array}{cccc} 
\xymatrix@=4mm{
&*+{1}\ar@{->}[ld]\\
*+{2}\ar@{->}[rr]&&*+{3}\ar@{<-}[lu]}  
&\xymatrix@=4mm{
&*+{1}\ar@{<-}[ld]\\
*+{2}\ar@{->}[rr]&&*+{3}\ar@{<-}[lu]}
&\xymatrix@=4mm{
&*+{1}\ar@{->}[ld]\\
*+{2}\ar@{<-}[rr]&&*+{3}\ar@{<-}[lu]}  
& \xymatrix@=4mm{
&*+{1}\ar@{->}[ld]\\
*+{2}\ar@{->}[rr]&&*+{3}\ar@{->}[lu]} 
\\
\xymatrix@=4mm{
&*+{1}\ar@{<-}[ld]\\
*+{2}\ar@{<-}[rr]&&*+{3}\ar@{<-}[lu]}  
& \xymatrix@=4mm{
&*+{1}\ar@{<-}[ld]\\
*+{2}\ar@{->}[rr]&&*+{3}\ar@{->}[lu]} 
&\xymatrix@=4mm{
&*+{1}\ar@{->}[ld]\\
*+{2}\ar@{<-}[rr]&&*+{3}\ar@{->}[lu]}  
&  \xymatrix@=4mm{
&*+{1}\ar@{<-}[ld]\\
*+{2}\ar@{<-}[rr]&&*+{3}\ar@{->}[lu]}  
\end{array}
\]
Since there are two choices of direction for each pair $(i, j)$ such that $1\le i<j\le n$, we have $|T_n|=2^{n(n-1)/2}$ in total.
}\eex






In what follows, the letter $G$ means an element of $T_n$ unless otherwise specified ($G$ is for Graph).
\bd{An \emph{inversion} of $G$ is a directed edge $j\to i$ with $j>i$.
The \emph{length} $\ell(G)$ is the number of inversions of $G$.
An \emph{upset} of an inversion $j\to i$ is a vertex $j$. 
Define $\om_G(j)$ to be the outdegree of $j$.
}\ed

\subsection{Cycle and transitivity}
Below, we just say a ``cycle" to mean a 3-cycle (which is the only kind of cycles we treat).


\bd{Let $(i, j, k)$ be a triple such that $i<j<k$.
Suppose $i, j, k$ form a cycle in $G$.
Say the cycle is \emph{positive} if $k\to j\to i\to k$; it is \emph{negative} if $i\to j\to k\to i$. Besides, say $G$ is \emph{transitive} if it does not contain any cycles.
}\ed

Observe that precisely six tournaments in Example \ref{t3} are transitive.

\subsection{Permutations}
By $S_n$ we mean the symmetric group on $[n]=\{1, 2, \dots, n\}$.
The set of \emph{inversions} of $w\in S_n$ is 
\[N(w)= \{(i, j)\in [n]\times [n]\mid \mbox{$i<j$ and $w^{-1}(i)>w^{-1}(j)$}\}.\]
Define the \emph{length} $\ell(w)$ to be $|N(w)|$.
Let $G(w)$ be the tournament such that $j\to i$ is an inversion of $G(w)$ $\iff (i, j)\in N(w)$.
Say the tournament $G(w)$ is \emph{induced} from a permutation $w\in S_n$.
Let us make sure the following:

\bft{Bressoud \cite[Exercise 2.4.2]{bressoud} There is a bijection between $S_n$ and transitive tournaments in $T_n$.
}\eft

Thanks to this result, we naturally view $S_n\sub T_n$ in what follows. In particular, $\ell(G(w))=\ell(w)$.





\subsection{Vandermonde determinant}

 \bd{Let $x_1, \dots, x_n$ and $\lam$ be commutative variables.
The \emph{$n$-th Vandermond $\lam$-determinant} is
\[V_n(\x, \lam)=
\prod_{1\le i<j\le n} (x_i+\lam x_j).\]
}\ed
This is a polynomial in $x_i$'s  ($\x$ means such variables for short) and $\lam$.
We must explain why we used the word ``determinant":
Following Robbins-Rumsey \cite{rr}, we recursively define a determinant-like function $|\phantom{A}|_{\lam}$ for square matrices as follows. First, we formally define $|\phantom{A}|_{\lam}$ for the $0$ by $0$ matrix to be $1$ and for a $1$ by $1$ matrix $(a_{11})$ to be $a_{11}$ itself. Now let $A$ be an $n$ by $n$ matrix for $n\ge 2$.
Let $A^j_i$ denote the matrix that remains when we delete the $i$-th row and $j$-th column of $A$.
If we wish to delete more than one row (column), the numbers of the deleted rows (columns) are listed as subscripts (superscripts).
The \emph{$\lam$-determinant} of $A$ is 
\[|A|_{\lam}=\frac{|A_1^1|_{\lam}|A_n^n|_{\lam}+\lam|A_n^1|_{\lam}|A_1^n|_{\lam}}{|A_{1n}^{1n}|_{\lam}} \mbox{  \quad(a rational function of $\lam$)}
\]
provided $|\phantom{A}|_{\lam}$ of all minors of $A$ are nonzero. 
In particular, $\lam=-1$ recovers the original determinant (going back to Dodgson and Desnanot-Jacobi).
From this point of view, we can understand $V_n(\x, \lam)$ as the $\lam$-determinant of the Vandermonde matrix: $V_n(\x, \lam)=|x_j^{n-i}|_{\lam}$.




\bd{
The \emph{Vandermonde monomial} for $G$ is $\rho(G)=\lam^{\ell(G)}\prod_{j\in [n]}x_j^{\om(j)}$.
}\ed

\bp{We have\label{sam}
\[V_n(\x, \lam)=\sum_{G\in T_n} \rho(G).\]
}\ep

\bpf{To a tournament $G$, assign a monomial with the choices of $x_i$ or 
$\lam x_j$ from each factor of $\prod_{i<j} (x_i+\lam x_j)$. Then $\lam$ in the monomial counts inversions and $x_j$ records the outdegree.

}\epf

Now, split the sum into two parts, transitive or not:
\[\sum_{G\in T_n} \rho(G)=\sum_{G\in S_n} \rho(G)+\sum_{G\in T_n\setminus S_n} \rho(G).\]

\bft{
\[\left.\sum_{G\in T_n\setminus S_n}\rho(G)\right|_{\lam=-1}=0.
\]
}\eft

\bpf{See Bressoud \cite[Exercise 2.4.4]{bressoud}.
}\epf

\section{Signed bigrassmannian statistic}
\label{too}

\subsection{Bigrassmannian statistic}

{\renewcommand{\arraystretch}{1.5}
\begin{table}[b]
\caption{statistics of $\ell$ and $\beta$ over $S_4$}
\begin{center}
\begin{tabular}{|c|c|c|c|c|c|c|c|c|c|c|c|c|c|c|c|}\hline
&$\ell$&$\beta$&&$\ell$&$\beta$&&$\ell$&$\beta$&&$\ell$&$\beta$\\\hline
1234&0&0&2134&1&1&3124&2&3&4123&3&6\\\hline
1243&1&1&2143&2&2&3142&3&5&4132&4&7\\\hline
1324&1&1&2314&2&3&3214&3&4&4213&4&7\\\hline
1342&2&3&2341&3&6&3241&4&7&4231&5&9\\\hline
1423&2&3&2413&3&5&3412&4&8&4312&5&9\\\hline
1432&3&4&2431&4&7&3421&5&9&4321&6&10\\\hline
\end{tabular}
\end{center}
\label{default}
\end{table}%
}

\bd{
Define the \emph{bigrassmannian statistic} for a tournament $G$ as 
\[
\beta(G)=\sum_{\substack{j\to i\\ j>i}}(j-i).\]
}\ed
Table \ref{default} shows this (and inversion) statistic over $S_4$.
\brk{
This statistic is named after the \emph{bigrassmannian permutations}; say $w\in S_n$ is \emph{bigrassmannian} if there exists a unique pair $(i, j)\in [n-1]\times [n-1]$ such that $w^{-1}(i)>w^{-1}(i+1)$ and $w(j)>w(j+1)$. We refer to Lascoux-Sch\"{u}tzenberger \cite{ls}, Geck-Kim \cite{geck}, Reading \cite{reading} and the author \cite{kob} for combinatorics of these permutations. 
}\erk

Define \emph{Bruhat order} $\le$ on $S_n$ as the transitive closure of the following binary relation: $v\to w$ meaning $w=vt_{ij}$, for some $i<j$, $t_{ij}$ a transposition and $\ell(v)<\ell(w)$.
Let $B(w)=\{u \mbox{ bigrassmannian}\mid u\le w\}$ and set $\beta(w)=|B(w)|$. The author \cite{kob} showed that  
\[\beta(w)=\sum_{(i, j)\in N(w)}(j-i).\]
Thus, we can compute $\beta$ simply as weighted enumeration of inversions:
\[\beta(3412)=(3-1)+(3-2)+(4-1)+(4-2)=8,\]
for example (Figure \ref{s33}). From this point of view, our definition above is a natural extension of $\beta$ for tournaments. 
This statistic implicitly appeared also in the Gessel-Viennot's lattice path counting context 
\cite[Theorem 3.7]{bressoud} as the quantity $\sum_{i=1}^ni(i-w(i))$:


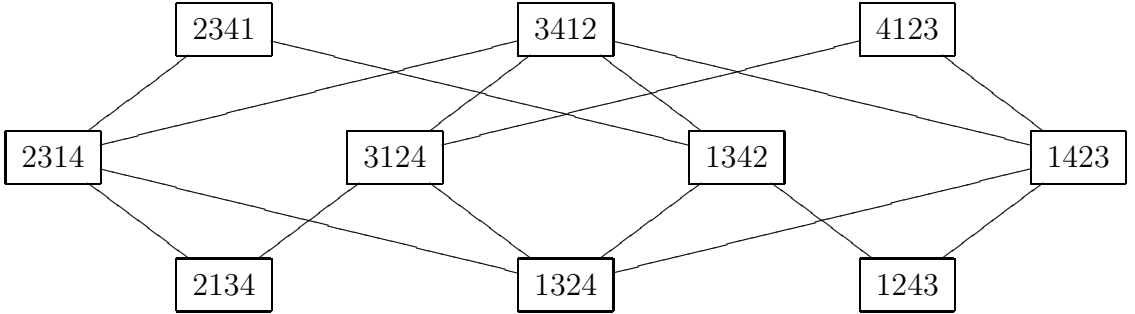
\begin{figure}[t]
\caption{Bruhat order of bigrassmannian permutation in $S_4$} \label{s33}

\[ \xymatrix@=10mm{
&*++<5pt>[F]{2341}\ar@{-}[dl]\ar@{-}[drrr]
&&*++<5pt>[F]{3412}\ar@{-}[dlll]\ar@{-}[dl]\ar@{-}[dr]\ar@{-}[drrr]
&&*++<5pt>[F]{4123}\ar@{-}[dlll]\ar@{-}[dr]&\\
*++<5pt>[F]{2314}\ar@{-}[dr]\ar@{-}[drrr]
&&*++<5pt>[F]{3124}\ar@{-}[dl]\ar@{-}[dr]
&&*++<5pt>[F]{1342}\ar@{-}[dl]\ar@{-}[dr]
&&*++<5pt>[F]{1423}\ar@{-}[dlll]\ar@{-}[dl]\\
&*++<5pt>[F]{2134}
&&*++<5pt>[F]{1324}
&&*++<5pt>[F]{1243}&\\
 }\]
\end{figure}


%

\bp{For each $w\in S_n$, we have
\label{bt}
\be{\item 
$\beta(w)=\sum_{i=1}^n (i-w(i))^2/2=\sum_{i=1}^n i(i-w(i))$.
\item $\beta(w)=\beta(w^{-1})$.
}\ee
}\ep

\bpf{(1) See \cite{kob} for the first equality.
It follows that \[\beta(w)=\fr{1}{\,2\,}\sum_{i=1}^n \left(i-w(i)\right)^2=\fr{1}{\,2\,}\sum_{i=1}^n\left(
i^2-2iw(i)+w(i)^2\right)
=\sum_{i=1}^n i(i-w(i)).\]
Next, (2) follows from the facts that (a) $u \mapsto u^{-1}$ is an order-preserving automorphism in Bruhat order on $S_n$, (b) $u$ is bigrassmannian $\iff$ so is $u^{-1}$; we do not go into details here because the proof is not so important for our discussions below.
}\epf
%







%




\bd{
Let $x_1, \dots, x_n$, $\lam$ and $q$ be commutative variables.
The \emph{weighted Vandermonde monomial} for $G$ 
is $\chi(G)=\lam^{\ell(G)}q^{\beta(w)}\prod_{j\in [n]}x_j^{\om(j)}$.
}\ed




\bd{The $n$-th \emph{weighted Vandermonde determinant} is 
\[V_n(\x, \lam, q)=\prod_{1\le i<j\le n}(x_i+\lam q^{j-i}x_j).
\]
}\ed
\renewcommand{\mu}{q}

\bex{
\label{y1}
\begin{align*}
V_3(x, \lam, \mu)&=(x_1+\lam \mu x_2)(x_1+\lam \mu^2x_3)(x_2+\lam \mu x_3)\\
&=x_1^2x_2+\lam \mu x_1x_2^2+\lam \mu x_1^2x_3+(\lam^2+\lam)\mu^2 x_1x_2x_3\\
&\phantom{=}+\lam^2 \mu^3 x_2^2x_3+\lam^2 \mu^3 x_1x_3^2+\lam^3 \mu^4 x_2x_3^2.
\end{align*}


}\eex

\bp{We have 
\[V_n(\x, \lam, q)=\sum_{G\in T_n}\chi(G).\]
}\ep

\bpf{
The idea is similar to Proposition \ref{sam}.
}\epf

\begin{lem}
\[\left.\sum_{G\in T_n\setminus S_n}\chi(G)\right|_{x_1=\cdots=x_n=1, \lam=-1}=0.
\]\label{this}
\end{lem}
To prove this lemma, we need a further definition and proposition.

\bd{For $i<j<k$, define a map $C_{ijk}: T_n\to T_n$
as follows:
if $i, j, k$ form a cycle in $G$, then 
$C_{ijk}(G)$ is the tournament with all three edges in the cycle reversed
and all other edges unchanged.
If $i, j, k$ do not form a cycle in $G$, then simply let $C_{ijk}(G)=G$.
}\ed
Observe that $C_{ijk}$ is an involution.


\bp{Let $i<j<k$.
If $i, j, k$ form a cycle in $G$, then $\ell(C_{ijk}(G))\in\{\ell(G)-1, \ell(G)+1\}$ and $\beta(C_{ijk}(G))=\beta(G)$.
}\ep

\bpf{A positive cycle contains two inversions whereas a negative cycle contains one.
The map $C_{ijk}$ interchanges these so that lengths differ by one. However, $\beta$ is invariant because of the equality $k-i=(k-j)+(j-i)$.
}\epf

\bpf[Proof of Lemma \ref{this}]{Consider the lexicographic order on $\{(i, j, k)\in [n]^3 \mid i<j<k\}$.
We will construct a perfect matching on the set $T_n\setminus S_n$.
First, choose all tournaments $G$ from $T_n$ such that $(1, 2, 3)$ is a cycle in $G$.
It is either positive or negative; hence $G \leftrightarrow C_{123}(G)$ gives a matching.
Next, choose all tournaments $H$ from the remaining tournaments such that $(1, 2, 4)$ is a cycle in $H$. Again, $H \leftrightarrow C_{124}(H)$ gives a matching.
Continue this procedure up to $(n-2, n-1, n)$.
We certainly exhausted all tournaments in $T_n\setminus S_n$ with the perfect matching constructed.
As shown above, each pair has lengths of opposite parity and the same $\beta$. 
Thus $x_1=\cdots =x_n=1$ and $\lam=-1$ yield zero.
}\epf


\subsection{Signed bigrassmannian polynomials}

\bd{Let $n$ be a positive integer. The $n$-th \emph{signed bigrassmannian polynomial} is 
\[B_n(q)=\sum_{w\in S_n}(-1)^{\ell(w)}q^{\beta(w)}.\]
}\ed

\bt{For all $n\ge 1$, we have 
\[B_n(q)=\prod_{k=1} ^{n-1}(1-q^k)^{n-k}.\]\label{mth1}
}\et

\bpf{

As before, split $V_n(\x, \lam, q)$ into two parts:
\[\prod_{1\le i<j\le n}(x_i+\lam q^{j-i}x_j)
=V_n(\x, \lam, q)=\sum_{G\in S_n}\chi(G) +\sum_{G\in T_n\setminus S_n}\chi(G).\]

With $x_1=\cdots =x_n=1$ and $\lam=-1$, the second sum vanishes as shown in Lemma \ref{this}.
As a result, we obtain
\[\prod_{1\le i<j\le n}(1-q^{j-i})=\sum_{w\in S_n}(-1)^{\ell(w)}q^{\beta(w)}\]
or 
\[B_n(q)=\prod_{k=1} ^{n-1}(1-q^k)^{n-k}.\]
}\epf

\begin{cor}
For $n\ge 3$, we have \[\sum_{w\in S_n}(-1)^{\ell(w)}\beta(w)=0.\]
In other words, the $\beta$-statistic is sign-balanced. 
\end{cor}

\bpf{
Note that $B_n(q)$ has a factor $(1-q)^{n-1}$ with $n-1\ge 2$.
Differentiate it once and let $q=1$. Then we get zero, as required.
}\epf

\bex{(cf. Reading's examples in Introduction)
\begin{align*}
B_2(q)&=1-q, \\
B_3(q)&=(1-q)^2(1-q^2)=1-2q+2q^3-q^4, \\
B_4(q)&=(1-q)^3(1-q^2)^2(1-q^3)\\
&=1-3q+q^2+4q^3-2q^4-2q^5-2q^6+4q^7+q^8-3q^9+q^{10}.
\end{align*}

}\eex

\section{Bigrassmannian determinant}
\label{sec4}
\subsection{Definition}
Next we want to understand $B_n(q)$ as a new sort of a determinant as mentioned in Introduction. From now on, we assume that $A=(a_{ij})=(a_{ij}(q))$ is an $n$ by $n$ matrix with entries 
being complex rational functions in $q^{1/2}$ (i.e., elements of $\cc(q^{1/2})$).
The reason why we introduce $q^{1/2}$ and $q^{-1}$ will be clearer in the next subsection.




\bd{The \emph{bigrassmannian determinant} of $A$ is  
\[\bdet(A)=\sum_{w\in S_n}(-1)^{\ell(w)}q^{\beta(w)}\prod_{i=1}^n a_{iw(i)}.\]
We formally define $\bdet$ of the $0$ by $0$ matrix to be $1$.
}\ed
For example, $\bdet(a_{11})=a_{11}$, 
$\bdet\left(\begin{array}{cc}a_{11} & a_{12} \\a_{21} & a_{22}\end{array}\right)=a_{11}a_{22}-qa_{12}a_{21}$ and
\begin{align*}\bdet\left(\begin{array}{ccc}a_{11}& a_{12}&a_{13} \\a_{21} & a_{22}&a_{23}
\\a_{31}& a_{32}&a_{33}
\end{array}\right)
&=a_{11}a_{22}a_{33}-qa_{12}a_{21}a_{33}-
q^{}a_{11}a_{23}a_{32}\\
&\phantom{=}+q^{3}a_{12}a_{23}a_{31}+
q^{3}a_{13}a_{21}a_{32}-
q^{4}a_{13}a_{22}a_{31}.
\end{align*}

\subsection{Matrix deformation}
We now give a more explicit description of the bigrassmannian determinant in terms of the original one. For this purpose, let us introduce a special term: a \emph{deformation} of $A=(a_{ij})$ is a new matrix $\bm{f}A:=(f_{ij}(q)a_{ij})$ for some indexed family of rational functions $\bm{f}=\{f_{ij}(q)\in \cc(q^{1/2})\mid (i, j)\in [n]\times [n]\}$. Note that the operation $a_{ij}\mapsto f_{ij}(q)a_{ij}$ may \emph{not} be $\cc(q^{1/2})$-linear in any rows nor columns. 
Hence it is in general difficult to predict how determinants change under such an operation.
However, as seen below, there are some nice cases:
%

\bd{Let $\bm{b}=\{b_{ij}(q)\}=\{q^{(i-j)^2/2}\}$. The \emph{bigrassmannian deformation} of $A$ is $\b A$.
}\ed

\bp{$\det(\b A)=\bdet(A)$.
}\ep

\bpf{By Proposition \ref{bt}, we have 
\begin{align*}
\det(\b A)&=\sum_{w\in S_n}(-1)^{\ell(w)}
\prod_{i=1}^nq^{(i-w(i))^2/2}a_{iw(i)}\\
&=\sum_{w\in S_n}(-1)^{\ell(w)}q^{\beta(w)}
\prod_{i=1}^na_{iw(i)}\\
&=\bdet(A).
\end{align*}
}\epf

\begin{thm}[a determinantal expression of $B_n(q)$] We have
\[B_n(q)=\det(q^{(i-j)^2/2}).\]\label{mth15}
\end{thm}

\bpf{$B_n(q)=\bdet(1)_{i, j=1}^n=\det(\b 1)=\det(q^{(i-j)^2/2})$.
}\epf
Observe determinantal expressions of $B_3(q)$ and $B_4(q)$:
\begin{align*}
\det\left(\begin{array}{ccc} 1 &q^{1/2}   & q^{4/2}  \\q^{1/2}  & 1  & q^{1/2}  \\q^{4/2}  &q^{1/2}   &1  \end{array}\right)&=(1-q)^2(1-q) \mbox{  and }\\
\det\left(\begin{array}{cccc} 1 & q^{1/2}
  & q^{4/2}
  & q^{9/2}
  \\q^{1/2}
  &  1 & q^{1/2}
  & q^{4/2}
  \\q^{4/2}
  &   q^{1/2}
&  1 & q^{1/2}
  \\ q^{9/2}
 & q^{4/2}
  &  q^{1/2}
 & 1 \end{array}\right)&=(1-q)^3(1-q^2)^2(1-q^3).
\end{align*}


We should now recognize that different deformations may give the same determinant:
given a family $\bm{f}$, there possibly exists $\bm{g}$ such that $\bm{g}\ne \bm{f}$ and 
$\det(\bm{f}A)=\det(\bm{g}A)$ for \emph{all} matrices $A$.
In particular, this is the case for $\b$:
Let $\b'=\{b_{ij}'\}=\{q^{i(i-j)}\}$ and $\b''=\{b_{ij}''\}=\{q^{j(j-i)}\}$ (here we need $q^{-1}$). 
Then $\det(\b A)=\det(\b'A)=\det(\b''A)$ as shown just below; since we could not find any references mentioning this little invariance, we here record it as a Corollary.

\begin{cor}
(little invariance of the determinant)
\[\det(q^{(i-j)^2/2}a_{ij})=\det(q^{i(i-j)}a_{ij})=\det(q^{j(j-i)}a_{ij}).\]
\end{cor}

\bpf{We only prove the first equality.
\begin{align*}
\det(q^{i(i-j)}a_{ij})&=\sum_{w\in S_n}(-1)^{\ell(w)}\prod_{i=1}^n
q^{i(i-w(i))}a_{i, w(i)}\\
&=\sum_{w\in S_n}(-1)^{\ell(w)}q^{\beta(w)}\prod_{i=1}^na_{i, w(i)}
=\det(q^{(i-j)^2/2}a_{ij}).
\end{align*}
}\epf
Such ``equivalent" deformations may be useful for evaluating and understanding combinatorial determinants (interpret $q^{(i-j)^2/2}$ as area of the triangle $(i, i)$, $(i, j)$ and $(j, j)$ in $\zz^2$); see Bressoud \cite[Section 3.3]{bressoud}, Gessel-Viennot \cite{gv} and Stembridge \cite{stembridge}, for details on Schur functions and nonintersecting lattice path counting by determinants. We will develop this idea in subsequent publications.


 

\subsection{Weighted condensation}

Our next task is to prove weighted condensation for bigrassmannian determinants; this is a natural idea as an analogy of the original determinant (and Robbins-Rumsey \cite{rr}).
Let $A$ be an $n$ by $n$ matrix with $n\ge 2$.
Recall that $A_{i}^j$ denotes the submatrix with the $i$-th row and $j$-th column deleted.
\bt{\label{mth2}
\[
\bdet(A)\bdet(A^{1n}_{1n})=
\bdet(A_1^1)\bdet(A_n^n)-q^{n-1}\bdet(A_n^1)\bdet(A_1^n).
\]
}\et

Some comments before the proof:
Let $A=(a_{ij}), C=\b A=(c_{ij})$ and $c_{ij}=q^{(i-j)^2/2}a_{ij}$.
For simplicity, we use $|\phantom{a_{ij}}|$ for the original determinant.

We will confirm the following five statements.
\be{\item $|C_{1n}^{1n}|=\bdet(A_{1n}^{1n})$.
\item $|C_1^1|=\bdet(A_1^1)$.
\item $|C_n^n|=\bdet(A_n^n)$.
\item $|C_n^1|=q^{(n-1)/2}\bdet(A_n^1)$.
\item $|C_1^n|=q^{(n-1)/2}\bdet(A_1^n)$.
}\ee 
Once we do this, then the conclusion follows from condensation for the original determinant:
\[|C||C^{1n}_{1n}|=
|C_1^1||C_n^n|-|C_n^1||C_1^n|.\]


\bpf{
(1) $|C_{1n}^{1n}|=\bdet(A_{1n}^{1n})$:
an $(i, j)$-entry of $C_{1n}^{1n}$ is $c_{i+1, j+1}$.
\begin{align*}
|C_{1n}^{1n}|&=|c_{i+1, j+1}|_{i, j=1}^{n-2}
=|q^{((i+1)-(j+1))^2/2}a_{i+1, j+1}|\\
&=
|q^{(i-j)^2/2}a_{i+1, j+1}|=\bdet (A_{1n}^{1n}).
\end{align*}
(2) $|C_1^1|=\bdet(A_1^1)$:
an $(i, j)$-entry of $C_{1}^1$ is $c_{i+1, j+1}$.
\begin{align*}
|C_1^1|&=|c_{i+1, j+1}|_{i, j=1}^{n-1}
=|q^{((i+1)-(j+1))^2/2}a_{i+1, j+1}|\\
&=|q^{(i-j)^2/2}a_{i+1, j+1}|=\bdet(A^1_1).
\end{align*}
(3) $|C_n^n|=\bdet(A_n^n)$:this is similar to (2).\\
(4) $|C_n^1|=q^{(n-1)/2}\bdet(A_n^1)$: an $(i, j)$-entry of $C_{n}^1$ is $c_{i, j+1}$.
\begin{align*}
|C_n^1|&=|c_{i, j+1}|_{i, j=1}^{n-1}=\det(q^{(i-(j+1))^2/2}a_{i, j+1})\\
&=|q^{((i-j)^2-2(i-j)+1)/2}a_{i, j+1}|
=q^{-\sum_{}i+\sum j}q^{(n-1)/2}|q^{(i-j)^2/2}a_{i, j+1}|=
q^{(n-1)/2}\bdet(A_n^1).
\end{align*}
(5) $|C_1^n|=q^{(n-1)/2}\bdet(A_1^n)$: this is similar to (4).

}\epf

Now we see an immediate consequence which is, however, not so obvious from the definition of $B_n(q)$. 
\begin{cor}\label{rec}
Signed bigrassmannian polynomials can be defined recursively as follows:
$B_1(q)=1, B_2(q)=1-q$ and $B_n(q)=\fr{\di {\strut B_{n-1}(q)^2}}{\di \strut B_{n-2}(q)} (1-q^{n-1})$
for $n\ge 3$.
\end{cor}

\bpf{
Apply the weighted condensation to $A=(1)_{i, j=1}^n$. All four determinants in the numerator are $B_{n-1}(q)$ while the denominator is $B_{n-2}(q)$.
}\epf

\section{Concluding remarks}
\label{sec5}

In this article, we introduced two new ideas, signed bigrassmannian polynomials and bigrassmannian determinant. 
We made use of tournaments as well as Vandermonde determinant to find $B_n(q)$. Then we introduced bdet as a $q$-analog of determinant as $q\to 1$ recovers the original one. Thanks to formulas of $\beta$-statistic, we obtained a determinantal expression of $B_n(q)$. Moreover, we established weighted condensation as an analogy of Robbins-Rumsey. 
After all, we did not find the unsigned statistic $\sum_{w\in S_n}q^{\beta(w)}$. Now an easy guess is to use the permanent instead. We leave this problem here for our future research.\\
We end with some more comments for subsequent work. 
\begin{itemize}
\item What is missing in our discussion is an alternating sign matrix (ASM) \cite{br3,robbins}.
Since inversions and bigrassmannian statistics also make sense for ASMs, we want to generalize some of our results to these matrices (note: we can extend $\beta$ for ASMs as the rank function of a distributive lattice). For example, what can we say about $\bdet$ for ASMs which are not permutations?.
\item We can also define ``$\lam q$-determinant" by replacing $\lam$ with $\lam q^{n-1}$ in Robbins-Rumsey condensation (provided all such minors are nonzero). Then we would obtain polynomials of the form $\prod (1+\lam q^{k})^{n-k}$, say $B_n(\lam, q)$. Then 
we can show as Corollary \ref{rec} that polynomials $\{B_n(\lam, q)\}$ satisfies 
\[B_n(\lam, q)=\di\fr{B_{n-1}(\lam, q)^2
}{B_{n-2}(\lam, q)}(1+\lam q^{n-1}).\]
Recently, there appeared such recursions and polynomials in the literature on Aztec diamonds, perfect matchings and domino tilings; see Brualdi-Kirkland \cite{brualdi}, Ciucu \cite{ciucu} and Elkies-Kuperberg-Larsen-Propp \cite{elkp}, for example. It would be nice to give an explicit connection between such work and our results.
\item As we mentioned Bruhat order, symmetric groups are Coxeter groups of type A. It makes sense to speak of a signed bigrassmannian statistic even in other situations: let $(W, S, \le)$ be a finite Coxeter system with Coxeter generators $S$ specified and $\le$ Bruhat order.  Define $\ell(w)=\min\{l\ge 0\mid w=s_1\cdots s_l, s_i\in S\}$ and the sign of $w$ to be $(-1)^{\ell(w)}$. Say $w$ is bigrassmannian if there exists a unique pair $(s_1, s_2)\in S\times S$ such that $\ell(s_1w)<\ell(w)$ and $\ell(ws_2)<\ell(w)$. Define $B(w)=\{u \mbox{ bigrassmannian} \mid u\le w\}$ and $\beta(w)=|B(w)|$ in the same way. Find a statistic $\sum_{w\in W} q^{\beta(w)}$.
\item We can think that each permutation $w$ gives a partition of an integer $\beta(w)$ with $\ell(w)$ parts as $\beta(w)=\sum_{(i, j)\in N(w)}(j-i)$; see Andrews-Eriksson \cite{andrews2} for the theory of integer partitions.
Then, it is natural to come to the following idea: \emph{Rothe diagram} for $w$ is the set $\{(i, j)\in [n]\times [n] \mid i < w^{-1}(j) \mbox{ and } j < w(i)\}$. As is well-known, the cardinality of this set is $\ell(w)$. Figure \ref{fult} shows an example; seven circles which does not cross any lines are elements of Rothe diagram for $w=35241$ (with $\beta(w)=15$). Is there any formula to compute $\beta$ from Rothe diagrams?

\begin{figure}[t]
\caption{Rothe diagram for 35241}
\label{fult}
\[\xymatrix@!=5mm{
*+{\circ}&*+{\circ}&*-{\bullet}\ar@{-}[rr]\ar@{-}[dddd]&*+{\circ}&*-{\circ}\\
*+{\circ}&*+{\circ}&*+{\circ}&*+{\circ}&*-{\bullet}\ar@{-}[ddd]\\
*+{\circ}&*-{\bullet}\ar@{-}[dd]\ar@{-}[rrr]&*+{\circ}&*+{\circ}&*-{\circ}\\
*+{\circ}&*+{\circ}&*+{\circ}&*-{\bullet}\ar@{-}[d]\ar@{-}[r]&*-{\circ}\\
*-{\bullet}\ar@{-}[rrrr]&*-{\circ}&*-{\circ}&*-{\circ}&*-{\circ}
}\]
\end{figure}
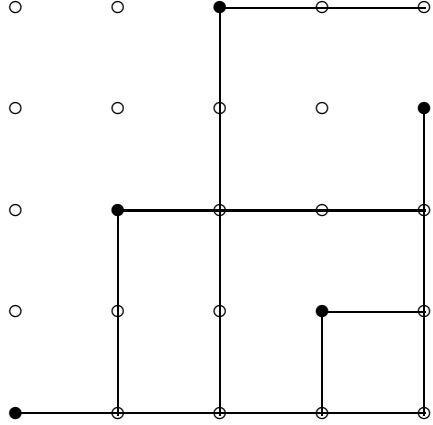



\end{itemize}



\section*{Acknowledgement}
The author thanks the anonymous referee for careful reading and advisory comments.
%
%
%
%
%
%

\end{document}